\documentclass{article}
\usepackage{amsmath}
\usepackage{theorem}
\usepackage{amssymb}
\usepackage{verbatim}
\usepackage{epsfig}
{\theorembodyfont{\rmfamily} \newtheorem{theorem}{Theorem}[section]

}
{ \theorembodyfont{\rmfamily} }
\numberwithin{equation}{section}
\begin{document}
\title{Note on the One and Two-Sided Z Tests}
\author{Terry R. McConnell\\
Syracuse University}
\maketitle
\begin{abstract}
The one sided Z test of elementary statistics is more powerful than the two-sided test of the same size.
\end{abstract}

\section{Introduction}

The humble Z test is taught in almost every introductory statistics course as a first illustration of hypothesis testing methodology. It assumes a normally distributed population with unknown mean $\mu$ and known variance $\sigma^2.$ The decision rule for the test rejects a null hypothesis of the form $H_0: \mu \le \mu_0$ in favor of the alternative $H_1: \mu > \mu_0$ if and only if $\bar x > \mu_0 + \frac{z_{\alpha}\sigma}{\sqrt{n}},$ where $\bar x$ is the mean of a random sample of size $n$, $\alpha$ is the desired size, and $\alpha \to z_{\alpha}$ is the inverse of the standard normal tail probability function. For the most common choice of $\alpha$ ($\alpha = 0.05$) one has $z_{\alpha} = 1.645,$ approximately.

The Z test is well known to be the best possible test of the given $H_0, H_1$ pair in the sense that its power function
$$
\pi_{\alpha}^+(\mu) = \Phi\left(\frac{\sqrt{n}}{\sigma}(\mu - \mu_0) - z_{\alpha}\right)
$$
exceeds that of any other test of the same alternatives in the region $\mu > \mu_0,$ assuming both tests have the same size. (See, for example, \cite{lehmann}.) The function $\Phi$ here is the cumulative distribution function of the standard normal distribution.

In the (symmetric) two-sided variant of the Z test one rejects
$H_0: \mu = \mu_0$ in favor of $H_1: \mu \ne \mu_0$ if and only if
$|\bar x - \mu_0| > z_{\frac{\alpha}2}\frac{\sigma}{\sqrt{n}}.$ The corresponding power function is given by
$$
\pi_{\alpha}^0(\mu) = 1 - \Phi\left(\frac{\sqrt{n}}{\sigma}(\mu - \mu_0) + z_{\frac{\alpha}2}\right) + \Phi\left(\frac{\sqrt{n}}{\sigma}(\mu - \mu_0) - z_{\frac{\alpha}2}\right).
$$
Many texts provide a graph similar to figure 1 below comparing the power functions of the one and 2-sided tests. (See, e.g., \cite[Figure 8.6]{degroot}.) Such graphs suggest that the one-sided test is always more powerful in the region corresponding to $H_1$, but this fact does not follow immediately from the optimality property of the one-sided test, since the one and 2-sided tests have different alternatives. Indeed, we have been unable to locate a proof in the literature. The purpose of this paper is to supply a proof that could be included in a calculus based introduction to mathematical statistics.

\begin{figure}
\begin{center}
\epsfig{file=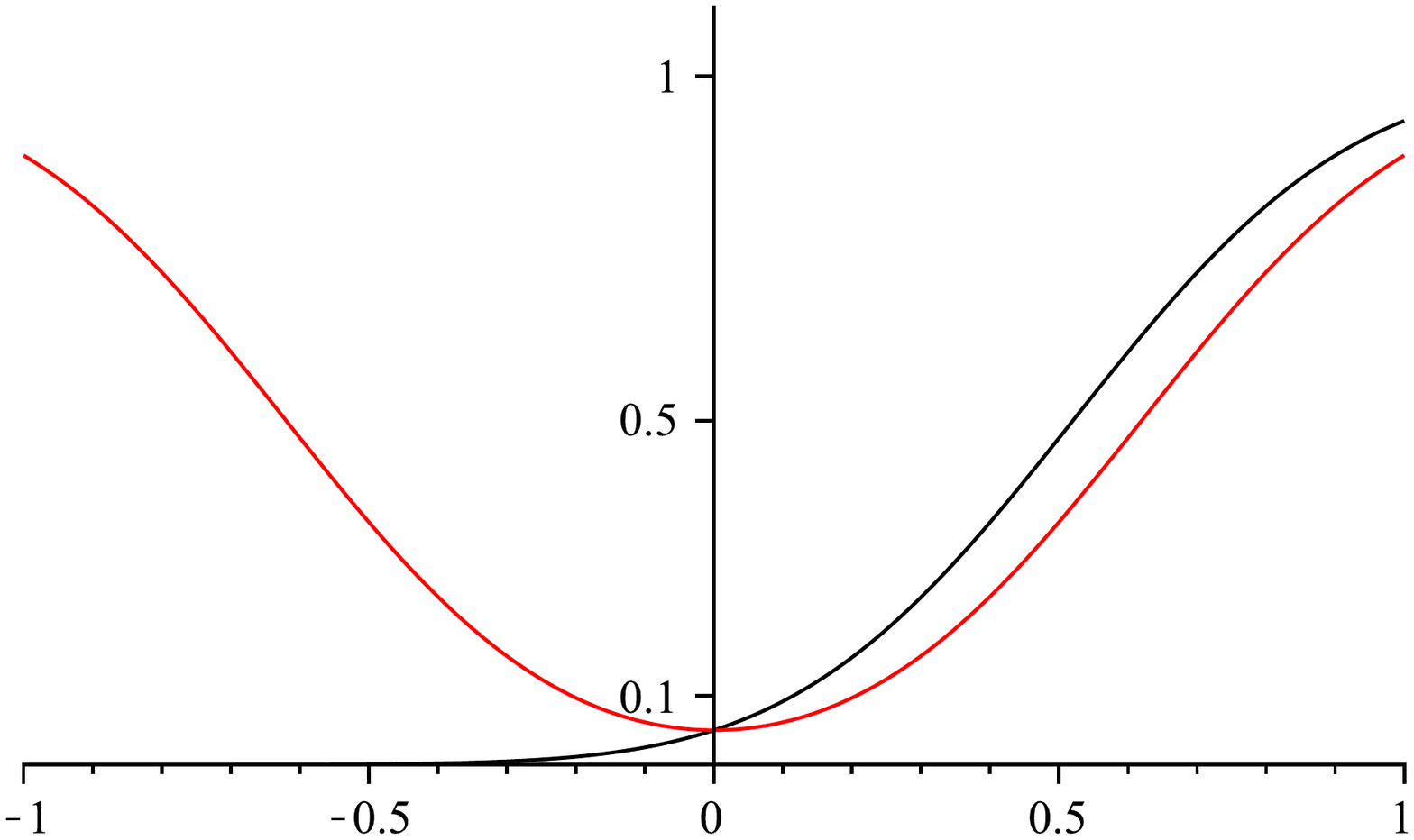,height=3in,width=3in}
\caption{One and 2-sided power functions with $n = 10, \alpha = 0.05,
\mu_0 = 0, $ and $\sigma = 1.$}
\end{center}
\end{figure}
In section 2 we prove that the one-sided Z test of $H_0: \mu \le \mu_0$ is more powerful than the 2-sided test of $H_0: \mu = \mu_0$, assuming both tests have the same size.
\begin{theorem} For each $0 < \alpha < 1$ we have
$$
\pi_{\alpha}^+(\mu) > \pi_{\alpha}^0(\mu),\ \ \mu > \mu_0.
$$
\end{theorem}
Section 3 is devoted to discussion of some related issues.

\section{Proofs}
Let
$$
\phi(x) = \frac1{\sqrt{2\pi}}e^{-\frac{x^2}2}
$$
be the standard normal probability density function, and $\Phi(x) = \int_{-\infty}^x \phi(t)\,dt$ the corresponding cumulative distribution function. After replacing $\frac{\alpha}2$ by $\alpha$, and $\frac{\sqrt{n}}{\sigma}(\mu - \mu_0)$ by $x$ for typographical convenience,  Theorem 1.1 reduces to the inequality
\begin{equation}
\Phi(x - z_{2\alpha}) > \Phi(x - z_{\alpha}) + \Phi(-x - z_{\alpha}), \ \ x > 0, 0 < \alpha < \frac12,
\end{equation}
for the function $\Phi.$

We shall provide 2 proofs of inequality (2.1). The first works for all values of $\alpha$ of practical interest, but not for all $\alpha.$ The second works for all $\alpha$.

For the first proof it is convenient to denote by $f(x)$ the function
$\Phi(x - z_{2\alpha}) - \Phi(x - z_{\alpha}),$ and by $g(x)$ the function $1 - \Phi(x + z_{\alpha}).$ In terms of these functions, inequality (2.1) becomes $f(x) > g(x), x > 0.$

Since $\phi$ is symmetric and unimodal, it is easy to see that the function $f$ is non-decreasing on the interval $(-\infty,\frac12(z_{2\alpha} + z_{\alpha})],$ and non-increasing on the interval $[\frac12(z_{2\alpha} + z_{\alpha}),\infty)$. Since $f(0) = \alpha = f(z_{2\alpha} + z_{\alpha})$, it follows that $f$ is bounded below by $\alpha$ for $0 \le x \le z_{2\alpha} + z_{\alpha}.$

On the other hand, $g$ is decreasing, and $g(0) = \alpha.$ Thus, it suffices to show that $f(x) > g(x)$ for $x > z_{2\alpha} + z_{\alpha},$ and we assume for the rest of the argument that $x$ lies in this range.
The Mean Value Theorem, applied on the interval $[x-z_{\alpha},x-z_{2\alpha}],$ gives the lower bound $f(x) > (z_{\alpha} - z_{2\alpha})\phi(x - z_{2\alpha}).$  On the other hand, the basic gaussian tail estimate gives the upper bound

$$
g(x) < \frac{\phi(x + z_{\alpha})}{z_{2\alpha}+2z_{\alpha}}.
$$

After a bit of algebra, one finds that the following inequality is sufficient:
\begin{equation*}
\frac{e^{\frac12(z_{2\alpha}^2-z_{\alpha}^2)}}{(z_{\alpha}-z_{2\alpha})
(z_{2\alpha}+2z_{\alpha})} < e^{x(z_{2\alpha}+z_{\alpha})}.
\end{equation*}
After replacing $x$ on the right hand side with its lower bound of $z_{2\alpha}+z_{\alpha}$ on the desired range, it is easy to check numerically that the resulting inequality holds for the ``classical'' values $z_{2\alpha} = 1.645, z_{\alpha} = 1.96,$ but fails to hold when $\alpha$ is sufficiently close to $\frac12$.

 Turning to the second proof of (2.1), since $\Phi(x - z_{2\alpha}) - \Phi(x - z_{\alpha})  = \int_{x - z_{\alpha}}^{x-z_{2\alpha}}\phi(t)\,dt$ and $\Phi(-x - z_{\alpha}) = \int_{x + z_{\alpha}}^{\infty}\phi(t)\,dt,$ it suffices to prove that
\begin{equation}
 \int_{x - z_{\alpha}}^{x-z_{2\alpha}}\phi(t)\,dt > \int_{x + z_{\alpha}}^{\infty}\phi(t)\,dt, \ \ x > 0.
\end{equation}
Let $\delta = z_{\alpha} - z_{2\alpha}.$ Then the left hand side of (2.2) minus the right hand side can be written as an integral over $(z_{2\alpha},z_{\alpha}]$:
\begin{equation}
\int_{z_{2\alpha}}^{z_{\alpha}}\phi(x - t) - \sum_{j=1}^{\infty}\phi(x + t + j\delta)\,dt.
\end{equation}
For $x > 0$ the integrand is equal to
\begin{equation*}
\phi(x)\left\{ e^{-\frac12t^2 + xt} - \sum_{j=1}^{\infty}e^{-\frac12(t+j\delta)^2}e^{-x(t+j\delta)}\right\}.
\end{equation*}
 The expression in brackets is clearly strictly increasing in $x$ for $t > 0$. Since $\phi(x) > 0$, the integral in (2.3) exceeds
\begin{equation*}
\phi(x)\int_{z_{2\alpha}}^{z_{\alpha}}e^{-\frac12t^2} - \sum_{j=1}^{\infty}e^{-\frac12(t+j\delta)^2}\,dt.
\end{equation*}
In turn, this expression equals
\begin{equation*}\sqrt{2\pi}\phi(x)\left\{\int_{z_{2\alpha}}^{z_{\alpha}}\phi(t)\,dt - \int_{z_{\alpha}}^{\infty}\phi(t)\,dt\right\},
\end{equation*}
which is equal to zero, by the definition of $z_{\alpha}$.
\section{Discussion}
The fact that one-sided tests are more powerful can lead to an ethical dilemma: A practitioner finds a 2 sided p-value of, say, 0.08, and `remembers' that they had intended to do a one-sided test all along. The new one-sided p-value of 0.04 is statistically significant.

Such flip-flops are not problematic, provided we insist that the practitioner provide a proof on {\it a priori} grounds that $\mu < \mu_0$ is impossible.

One might wonder whether the result of Theorem 1.1 is a peculiar feature of the normal distribution, or whether there is some deeper phenomenon at work that is applicable to more general populations? It is easy to find distribution families for which there appear to be no obvious analogues of Theorem 1.1. Consider, for example, the translated Cauchy distributions. Let $Y$ be a random variable having probability density function given by
$$
h(x) = \frac1{\pi} \frac1{1 + x^2}, -\infty < x < \infty,
$$
and let $X = m + Y$, where the real number $m$ serves as the parameter of interest. We may treat the random variable $X$ as a sample of size $n = 1.$ To test $H_0: m \le 0$ versus $H_1: m > 0,$ the analogue of the Z test rejects $H_0$ if and only if $X > c_{\alpha}$, where $c_{\alpha} = \cot(\pi\alpha).$ Similarly, the 2-sided test of $H_0: m = 0,$ versus $H_1: m \ne 0$ rejects $H_0$ if and only if $|X| > c_{\frac{\alpha}2}.$ The one-sided test has power function
$$
\frac12 + \frac1{\pi}\tan^{-1}(m - c_{\alpha}),
$$
and the 2-sided test has power function
$$
1 - \frac{\tan^{-1}\left(m + c_{\frac{\alpha}2}\right) -
\tan^{-1}\left(m - c_{\frac{\alpha}2}\right)}{\pi},
$$
so a result analogous to Theorem 1.1 would entail
$$
\tan^{-1}\left(m - c_{\alpha}\right) + \tan^{-1}\left(m + c_{\frac{\alpha}2}\right) > \frac{\pi}2 + \tan^{-1}\left(m - c_{\frac{\alpha}2}\right), m > 0, 0 < \alpha < 1.
$$
This inequality, however, is not always true. For example, it fails when $\alpha = \frac12$ and $ m = 2.$ (In that case, equality holds.)

\section{Acknowledgement} It is a pleasure to thank Hyune-Ju Kim for a conversation related to the first proof in section 2.

\end{document}